\documentclass[11pt, letterpaper]{article}

\usepackage[utf8]{inputenc}

\usepackage[
letterpaper,
top=1in,
bottom=1in,
left=1in,
right=1in]{geometry}

\usepackage{amsthm, amsmath, amssymb}
\usepackage{mathtools}
\usepackage{xcolor}
\usepackage{float}

\usepackage[colorlinks=true,urlcolor=blue,linkcolor=blue,citecolor=blue]{hyperref}
\usepackage[nameinlink]{cleveref}
\crefname{equation}{}{}

\usepackage{tikz}
\usepackage{graphicx}

\newtheorem{theorem}{Theorem}[section]
\newtheorem{corollary}[theorem]{Corollary}
\newtheorem{lemma}[theorem]{Lemma}

\newtheorem{assumption}[theorem]{Assumption}
\newtheorem{conjecture}[theorem]{Conjecture}
\newtheorem{proposition}[theorem]{Proposition}

\theoremstyle{definition}
\newtheorem{definition}[theorem]{Definition}
\newtheorem{example}[theorem]{Example}

\newtheorem{remark}[theorem]{Remark}

\newcommand{\R}{\mathbb{R}}
\newcommand{\N}{\mathbb{N}}

\newcommand{\mC}{\mathcal{C}}

\newcommand{\CDco}[2]{\Lambda_{\otimes #1}^{#2}}
\newcommand{\CDcoreg}[2]{\widetilde\Lambda_{\otimes n}^{#2}}
\newcommand{\CDtot}[2]{\Lambda_{#1}^{#2}}
\newcommand{\CDtotreg}[2]{\widetilde\Lambda_{#1}^{#2}}
\newcommand{\sparseCD}[3]{\Gamma_{#1}^{#2, #3}}
\newcommand{\sparseCDshort}[3]{\Gamma_{#1}^{#2}}
\newcommand{\sparseCDreg}[3]{\widetilde\Gamma_{#1}^{#2, #3}}
\newcommand{\sparseCDregshort}[3]{\widetilde\Gamma_{#1}^{#2}}

\newcommand{\Rco}{\mathbb{R}[\x]_{\otimes n}}
\newcommand{\tens}[1]{\otimes{#1}}

\newcommand{\mainset}{\mathbf{\Omega}}
\renewcommand{\epsilon}{\varepsilon}

\newcommand{\x}{\mathbf{x}}
\newcommand{\y}{\mathbf{y}}
\newcommand{\z}{\mathbf{z}}

\newcommand{\dist}{\mathrm{dist}}
\newcommand{\clo}[1]{\mathrm{cl}(#1)}
\newcommand{\inter}[1]{\mathrm{int}(#1)}
\newcommand{\clique}{\mathrm{clq}}

\title{A Christoffel-like function for high-dimensional support inference in graphical models}

\author{
  \begin{minipage}[t]{0.5\textwidth}
    \centering
    Jean B. Lasserre\thanks{The first author is supported by the AI Interdisciplinary Institute ANITI  funding through the French program
``Investing for the Future PI3A" under the grant agreement number ANR-19-PI3A-0004. This research is also part of the programme DesCartes and is supported by the National Research Foundation, Prime Minister's Office, Singapore under its Campus for Research Excellence and Technological Enterprise (CREATE) programme.} \\ \small \href{mailto:lasserre@laas.fr}{lasserre@laas.fr} \\ \small LAAS-CNRS \& University of Toulouse
  \end{minipage}
  \hfill
  \begin{minipage}[t]{0.45\textwidth}
    \centering
    Lucas Slot\thanks{This research was completed while the second author was at ETH Zurich, Institute of Theoretical Computer Science. 
    }  
     \\ \small \href{mailto:l.f.h.slot@uva.nl}{l.f.h.slot@uva.nl} \\ \small University of Amsterdam
  \end{minipage}
}

\begin{document}

\maketitle
            
\begin{abstract}
Christoffel polynomials are classical tools from approximation theory. They can be used to estimate the (compact) support of a measure $\mu$ on~$\R^d$ based on its low-degree moments. Recently, they have been applied to problems in data science, including outlier detection and support inference. A major downside of Christoffel polynomials in such applications is the fact that, in order to compute their coefficients, one must invert a moment matrix whose size grows rapidly with the dimension $d$. In this paper, we propose a modification of the Christoffel polynomial which is significantly cheaper to compute, but retains many of its desirable properties. In particular, it (1) exhibits a so-called support dichotomy and (2) it is a rational function, whose numerator and denominator factor into `lower-dimensional' Christoffel polynomials whose coefficients can be computed by inverting potentially much smaller moment matrices.
Our approach relies on sparsity of the underlying measure $\mu$, described by a graphical model. The complexity of our modification depends on the treewidth of this model.
\end{abstract}

\section{Introduction} \label{SEC:INTRO}
The Christoffel-Darboux (CD) kernel and the closely related Christoffel polynomial
are classical tools from the theory of approximation and orthogonal polynomials.
One of their main applications is the study of a measure $\mu$ with compact support $\mainset \subseteq \R^d$ using only knowledge of its moments. 
For each $n \in \N$, we write~$\CDtot{n}{\mu}$ for the Christoffel polynomial of degree $2n$ associated to $\mu$.
Writing $P_\alpha$, $\alpha \in \N^d$ for an orthonormal system of polynomials w.r.t. $\mu$ with $\mathrm{deg}(P_\alpha) = \|\alpha\|_1$, it equals\footnote{The Christoffel \emph{polynomial} is the diagonal of the Christoffel-Darboux \emph{kernel}, which is given by $\mathrm{CD}_n(\x, \y) = \sum_{\|\alpha\|_1 \leq n} P_\alpha(\x) P_\alpha(\y)$.
It is the reciprocal of the Christoffel \emph{function} which is also found in the literature.}
\begin{equation} \label{EQ:Cpol-intro}
    \CDtot{n}{\mu}(\x) = \sum_{\|\alpha\|_1 \leq n} P_\alpha(\x)^2\,, \quad \forall \x \in \R^d.
\end{equation}
See \Cref{SEC:CDprelim} for a detailed definition (and~\Cref{SEC:notations} for an overview of notations). 
As described there, and in~\Cref{SEC:CDsupdich}, Christoffel polynomials have the following features:
\begin{enumerate}
    \item The polynomial $\CDtot{n}{\mu}$ is a \emph{sum of squares} of polynomials of degree~$n$;
    \item The coefficients of $\CDtot{n}{\mu}$ can be computed by inverting a matrix of size ${n + d \choose d}$, whose entries are the \emph{moments} of $\mu$ of degree  at most $2n$;
    \item For $\x \in \R^d$, the polynomials $\CDtot{n}{\mu}$ exhibit a \emph{support dichotomy}:
\[
    \CDtot{n}{\mu}(\x) \sim \begin{cases} \exp(\alpha \cdot n), ~ \alpha>0 \quad & \text{when } \x \not\in \mainset,\\ 
    \mathrm{poly}(n) & \text{when } \x \in \mathrm{int}(\mainset), \end{cases} \quad \text{ as $n \to \infty$.}
\]
\end{enumerate}
For sufficiently large $n$, support dichotomy tells us that $\CDtot{n}{\mu}(\x)$ is relatively big for $\x \not\in \mainset$, but relatively small for $\x \in \mathrm{int}(\mainset)$. 
The \emph{sublevel sets} of $\CDtot{n}{\mu}$ may therefore be used to identify the support $\mainset$ of~$\mu$, which leads to \emph{semialgebraic} approximations of~$\mainset$. 
It was only recently \cite{LasserrePauwels2019,LPP:2022} that this simple property proved to be very useful to help solve some important problems in data mining and data analysis. Since then, the CD kernel and Christoffel polynomial have found applications in support~estimation~\cite{Devonport:supportestimation}, outlier detection~\cite{kernalizedCD, Ducharlet:outlier}, stochastic differential equations~\cite{Henrion:SDE}, functional approximation~\cite{Marx:CDapproximation} and
topological data analysis~\cite{RoosHoefgeestSlot:TDA}.
Quoting~\cite{LPP:2022}:
\begin{quote}
\emph{Among the many positive definite kernels appearing in classical analysis, approximation theory, probability, mathematical physics, control theory and more recently in machine learning, Christoffel-Darboux kernel (CD kernel in short) stands aside by its numerical accessibility from raw data and its versatility in encoding/decoding fine properties of the generating measure.}
\end{quote}

\subsection*{Computational (in)tractability}
The primary bottleneck for practical computations using Christoffel polynomials is the rapidly increasing size of the moment matrix that needs to be inverted to find their coefficients. 
Indeed, this matrix is of size ${n + d\choose d}$, which is prohibitive already for moderate values of $n, d$.
The goal of this paper is to modify the Christoffel polynomial in such a way that it becomes significantly easier to compute when the dimension $d$ is large, while retaining its desirable properties. As we will see, so far we are able to make such modifications only in the case where the measure $\mu$ has some additional structure.

\subsection*{Product measures}
The motivating example is when 
$\mu = \mu_1 \otimes \mu_2$ is a \emph{product measure} on~${\R^{d_1} \times \R^{d_2}}$.
Surprisingly, it is not clear whether a more efficient strategy to compute the Christoffel polynomial $\CDtot{n}{\mu}$ exists in this setting.
In~\Cref{DEF:modifiedCD} below, we propose a straightforward modification $\CDco{n}{\mu}$ of $\CDtot{n}{\mu}$ which, in this setting, decomposes into the product
\begin{equation} \label{EQ:productdecomp-1}
    \CDco{n}{\mu}(\x) = \CDco{n}{\mu_1}(\x_1)\cdot\CDco{n}{\mu_2}(\x_2)
    \,,\quad\forall \x = (\x_1, \x_2) \in\R^d\,;\:\forall n\in\N.
\end{equation}
The upshot is that, in order to compute $\CDco{n}{\mu}$ it now suffices to invert moment matrices associated to the (lower dimensional) measures $\mu_1$ and $\mu_2$.
When ${d=1}$, our modification $\CDco{n}{\mu}$ agrees with the standard Christoffel polynomial $\CDtot{n}{\mu}$. Therefore, if $\mu=\otimes_{i=1}^d \mu_i$ is a product of univariate measures, it decomposes into 
a product of univariate (standard) Christoffel polynomials:
\begin{equation} \label{EQ:productdecomp}
    \CDco{n}{\mu}(\x) = \prod_{i=1}^d \CDtot{n}{\mu_i}(x_i)\,,\quad\forall \x\in\R^d\,;\:\forall n\in\N.
\end{equation}
The coefficients of each factor in the decomposition \eqref{EQ:productdecomp} can be computed by inverting a matrix whose size \emph{does not depend on $d$} at all.

The key difference between our modification $\CDco{n}{\mu}$ and~$\CDtot{n}{\mu}$ (which makes the above decompositions possible) is that the former is given by a sum of the type~\eqref{EQ:Cpol-intro} ranging over $\alpha$ with $\|\alpha\|_\infty \leq n$ (rather than $\|\alpha\|_1 \leq n$). As a result, it is a polynomial of \emph{coordinate-wise} (rather than \emph{total}) degree $2n$, i.e., using monomials $\x^\alpha$, $\alpha \in \N^n$ with $\|\alpha\|_\infty \leq 2n$. Subsequently, computing the coefficients of $\CDco{n}{\mu}$ naively requires inverting a moment matrix of size $(n+1)^d$, versus ${n+d \choose d}$ for $\CDtot{n}{\mu}$. Still, in light of~\eqref{EQ:productdecomp-1} and~\eqref{EQ:productdecomp} computing $\CDco{n}{\mu}$ could be far cheaper when $\mu$ has product structure. 
Importantly, the \emph{coordinate-wise Christoffel polynomial} $\CDco{n}{\mu}$ is still a sum of squares, and it is not hard to show that it retains the support dichotomy property (see \Cref{SEC:modifiedCD:asymptotics}).

\subsection{Main contribution: a rational Christoffel-like function}
The assumption that $\mu$ is a product measure is very restrictive.
The primary contribution of this paper is to extend the decomposition~\eqref{EQ:productdecomp} to the more general class of measures that exhibit a \emph{conditional} product structure represented by a \emph{graphical model} $G$ (see \Cref{SEC:graphmodels}). 
The (computational) complexity of our approach will depend on the \emph{treewidth} $\tau(G)$ of this model, which captures in some sense how close $\mu$ is to being a product measure. (For instance, if $\mu$ is a product measure, then $\tau(G) = 0$.)
Concretely, given a measure $\mu$ with graphical model $G$, we propose in~\Cref{DEF:rational-cd} below a sequence of \emph{rational functions} $\sparseCDshort{n}{\mu}{G, J}$, $n \in \N$, of the form
\begin{equation}
    \label{EQ:rational-intro}
\sparseCDshort{n}{\mu}{G, J}(\x) = \sparseCD{n}{\mu}{G, J}(\x) := \frac{\prod_{v \in V}\CDco{n}{\mu_{A_v}}(\x_{A_v})}{\prod_{e \in E}\CDco{n}{\mu_{B_e}}(\x_{B_e})}\,,\quad \x\in\R^d.
\end{equation}
Here, the measures $\mu_{A_v}, \mu_{B_e}$, $v \in V$, $e \in E$ are the \emph{marginals} of $\mu$ on the variables indexed by (possibly overlapping) subsets $A_v, B_e \subseteq [d]$ that correspond respectively to the vertices and edges of a so-called \emph{junction tree} $J = (V, E)$ associated with~$G$. The vertices of this tree represent maximal cliques in a \emph{chordal completion} of $G$, and its edges represent intersections between these cliques. See~\Cref{SEC:graphmodels} for the precise definition. For now, let us emphasize that $\sparseCDshort{n}{\mu}{G, J}$ depends on the given graphical model~$G$, and on our choice of junction tree for $G$, which is not unique in general.
Our \emph{rational Christoffel-like functions}\footnote{This terminology refers to the fact that $\sparseCDshort{n}{\mu}{G, J}$ is a rational function, whose numerator and denominator are (products of) Christoffel polynomials. It is not to be confused with the \emph{Christoffel function} found commonly in the literature (being the reciprocal of the Christoffel polynomial).} have the following distinguishing properties:
\begin{enumerate}
    \item The factors of $\sparseCDshort{n}{\mu}{G, J}$ are (modified) Christoffel polynomials associated with marginals of $\mu$. In particular, they are sums of squares of polynomials of coordinate-wise degree $n$;
    \item The coefficients of each factor of $\sparseCDshort{n}{\mu}{G, J}$ can be computed by inverting a matrix of size at most $(n+1)^{\clique(J)}$, whose entries are moments of $\mu$. 
    Here, $\clique(J) = \max_{v \in V} |A_v|$ denotes the size of a largest clique represented by the junction tree $J$. 
    If $J$ is chosen optimally (i.e., corresponding to an optimal chordal completion of $G$), then $\clique(J) = \tau(G) + 1$. In total, there are at most~$2d$ such factors present in~\eqref{EQ:rational-intro}.
    \item Under mild assumptions on $\mu$ (see \Cref{ASSU:measure}), for almost all $\x \in \R^d$, the functions $\sparseCDshort{n}{\mu}{G, J}$ exhibit a support dichotomy:
    \[
    \sparseCDshort{n}{\mu}{G, J}(\x) \sim \begin{cases} \exp(\alpha \cdot n), ~ \alpha > 0 \quad & \text{when } \x \not\in \mainset,\\ 
    \mathrm{poly}(n) & \text{when } \x \in \mathrm{int}(\mainset), \end{cases} \quad \text{ as $n \to \infty$.}
    \]
\end{enumerate}
So, while on the one hand our rational Christoffel functions are significantly easier to compute than standard Christoffel polynomials when $\tau(G) \ll d$, on the other hand, they also retain their desirable support dichotomy property, and their sublevel sets are still semialgebraic. The dichotomy property is crucial for support inference in data analysis and mining. For instance the Christoffel polynomial 
has been used as a simple score function to detect outliers in e.g. \cite{Ducharlet:outlier,LasserrePauwels2019}. In the present context one may also use $\sparseCDshort{n}{\mu}{G,J}$ in \eqref{EQ:rational-intro} as a score function to detect outliers without computing 
the costly Christoffel polynomial of the joint distribution.
We establish support dichotomy for $\sparseCDshort{n}{\mu}{G,J}$ in \Cref{SEC:ratfunction:supdich}. This is our main technical result~\Cref{THM:sparsesupdich}. For its proof, we rely on the \emph{clique-intersection property} of junction trees (also known as \emph{running intersection property}).

Finally, we note that $\sparseCDshort{n}{\mu}{G, J}$ is a natural generalization of the modified Christoffel polynomial $\CDco{n}{\mu}$ in the following sense: if $\mu = \otimes_{i=1}^d \mu_i$ is a product measure, then, \emph{regardless of the provided graphical model for $\mu$ and our choice of junction tree}, we have
\[
\sparseCDshort{n}{\mu}{G, J}(\x) = \CDco{n}{\mu}(\x) = \prod_{i=1}^d \CDtot{n}{\mu_i}(x_i)\,,\quad\forall \x\in\R^d\,;\:\forall n\in\N.
\]
See \Cref{PROP:agreewhenproduct} for details. Intuitively, this means that our rational Christoffel function always recovers the (modified) Christoffel polynomial when $\mu$ is a product measure, even if that information was not reflected by the graphical model.

\subsection{Second contribution: density estimation}
If~$\mu$ has a density $f$ w.r.t. the Lebesgue measure on~$\mainset$, Christoffel polynomials can be used to approximate $f$. Namely, under certain regularity conditions on $\mainset$ and $\mu$, it can be shown that
\[
    \lim_{n \to \infty} {n + d \choose d}^{-1} \cdot \, \CDtot{n}{\mu}(\x) = \frac{\omega_E(\x)}{f(\x)},
\]
where $w_E$ is the so-called \emph{equilibrium measure} of $\mainset$ (which is independent of~$f$). 
This approximation is useful only if $\omega_E$ is known a priori, which is not the case in general. 
To circumvent this issue, Lasserre~\cite{Lasserre:regularizedCD} introduced a \emph{regularization}~$\CDtotreg{n}{\mu}$ of the Christoffel polynomial, which depends on an additional variable~${\epsilon > 0}$ (see \Cref{SEC:cdregular}). 
Under a (Lipschitz) continuity assumption on $f$, for any $\x \in \mathrm{int}(\mainset)$, it satisfies
\[
    \lim_{n \to \infty} \epsilon^{-d} \cdot \CDtotreg{n}{\mu}(\x; \epsilon)^{-1} = {f(\x) + O(\epsilon)}, \quad \text{ as } \epsilon \to 0.
\]See, \cite[Theorem 2.5]{Lasserre:regularizedCD}. Here, the $O$-notation hides constants that depend on $\x$ and $f$.
It thus approximately recovers $f(\x)$ when $\epsilon$ is small (without any unknown factors). 
In \Cref{SEC:ratfunction:regular}, we define an analogous regularization $\sparseCDregshort{n}{\mu}{G, J} = \sparseCDreg{n}{\mu}{G, J}$ of our rational Christoffel functions. We show that this allows for a similar recovery of the density:
\[
    \lim_{n \to \infty} \epsilon^{-d} \cdot \sparseCDregshort{n}{\mu}{G, J}(\x;\varepsilon)^{-1} = {f(\x) + O(\epsilon)}, \quad \text{ as } \epsilon \to 0.
\]
See \Cref{THM:sparseconvdensity} and \Cref{COR:sparsedensityeps}. As before, the advantage is that computation of~$\sparseCDregshort{n}{\mu}{G, J}$ requires inversion of far smaller moment matrices than $\CDtotreg{n}{\mu}$ when $\tau(G) \ll d$.

\subsection{Related work}

\medskip
\emph{Graphical models.}
Graphical models (or Bayesian Networks) are a classical tool in computational statistics, see, e.g.,~\cite{ProbabilisticGraphicalModels:book, jordan} 
for a general reference. A typical application for probabilistic inference
is efficient computation of marginal probabilities of a known high-dimensional joint distribution, a NP-hard problem in general as discussed in e.g. \cite[\S 1]{marvin}.
For example, given a (discrete) distribution $D$ on $\{-1, 1\}^d$, computing $\ {P}_{\x \sim D}(x_1 = 1)$ naively requires summation of $2^{d-1}$ terms:
\[
    \mathbb{P}_{\x \sim D}(x_1 = 1) \quad = \sum_{y_2, \ldots, y_d \in \{-1, 1\}} \mathbb{P}_{\x \sim D}(x_1=1, x_2=y_2, \ldots, x_d=y_d).
\]
However, when $D$ satisfies some appropriate conditional independence relations (expressed by a graphical model), much more efficient algorithms exist to compute or approximate such probabilities; see e.g. 
the use of quadrature rules in \cite{marvin} and Kikuchi approximations in \cite{kikuchi-2}. Learning parameters when data contain missing values is also an important research issue in Bayesian Networks (see e..g Bethe/Kikuchi approximation methods described in \cite{kikuchi-1}).

Morally speaking, this work differs from the above type of application in two ways. First, in our setting, the underlying measure $\mu$ is unknown, and only its moments are available. Second, computing the Christoffel polynomial~$\CDtot{n}{\mu_A}$ for a marginal $\mu_A$, $A \subseteq [d]$ of $\mu$ is actually \emph{easier} than computing~$\CDtot{n}{\mu}$. 
Indeed, the moment matrix $M^{\mu_A}_n$ of $\mu_A$ is a submatrix of~$M^{\mu}_n$. For us, the goal is in fact to compute, approximate or emulate $\CDtot{n}{\mu}$ using products and ratios of  Christoffel polynomials associated to marginals $\mu_A$ of $\mu$ with $|A| \ll d$. For instance, this strategy allows to detect outliers with the same methodology as in \cite{LasserrePauwels2019}, but without the need to obtain the (costly) inverse of the whole moment matrix~$M^\mu_n$.

\medskip
\emph{Sparse approximation of Christoffel polynomials.}
In~\cite{kernalizedCD}, the authors propose a method to (approximately) evaluate the Christoffel polynomial for large $d$ under a type of sparsity assumption on $\mu$. Namely, they assume that $\mu$ is a discrete measure supported on~$N$ atoms $\y_1, \ldots \y_N \in \R^d$. For instance, $\mu$ could be the \emph{empirical measure} associated to a draw of $N$ independent samples according to some probability measure on $\R^d$. 
The upshot is that, in this case, the moment matrix of $\mu$ has rank at most $N$. For fixed $\x \in \R^d$, and $\rho > 0$ this allows computation of a lower bound $\varphi(\x; \rho) \leq \rho \CDtot{n}{\mu}(\x)$ on the Christoffel polynomial in time $O(N^3 + N^2d\log n)$, see~\cite{kernalizedCD} for details. 
We note that the conditional independence assumption we make on $\mu$ in the present work (expressed by a graphical model) is not comparable to the sparsity assumption made in~\cite{kernalizedCD}. In particular, the sparsity in \cite{kernalizedCD} is related to a small discrete support in $\R^d$ while ours is related to a certain number of 
marginal measures $\mu_i$ on spaces $\R^{d_i}$ (each of smaller dimension $d_i\ll d$) associated with cliques in some graph
representing interactions between variables.
Moreover, it applies to non-discrete measures.

\medskip \emph{Sparse polynomial optimization.} In polynomial optimization, one is tasked with finding the minimum value attained by a polynomial \emph{cost function}, on a \emph{feasible region} in $\R^d$ described by polynomial inequalities. In large-scale polynomial optimization, it is often the case that not all variables interact directly in the description of the problem. Rather, they are weakly coupled: 
\begin{itemize}
    \item the cost function is a sum of 
    polynomials that each depend  only on a few variables;
    \item the feasible region is described by potentially many polynomial constraints each of which also involves only a few variables. 
\end{itemize}
    Such interactions between the variables can be represented by a graphical model and its associated junction tree analogously to conditional product structure of measures (as described in \Cref{SEC:graphmodels} below). Using this representation, 
    the authors in \cite{Waki2006} have defined a \emph{sparsity-adapted} version of the celebrated \emph{Moment-SOS hierarchy}, which (in the non-sparse setting) provides a converging sequence of \emph{lower bounds} on the solution of a polynomial optimization problem.
    Their sparse hierarchy can handle optimization problems with a large number of variables provided that  
    the graphical model has small treewidth. Moreover and importantly, the resulting bounds still
    converge towards the global optimum. The rationale behind the convergence in the sparse setting is a sparsity-adapted certificate of positivity \emph{\`a la Putinar} provided in \cite{Lasserre06a}. It relies critically on the clique-intersection property of junction trees (referred to as \emph{running intersection property} in~\cite{Lasserre06a}). This mirrors our reliance on the clique-intersection property to show support dichotomy for the rational Christoffel functions in~\Cref{SEC:ratfunction:supdich}. 
    
    For further reading on sparsity in polynomial optimization, see, e.g.,~\cite{HKLbook,Magron-sparse}.

\subsection*{Outline} The rest of the paper is organized as follows. In \Cref{SEC:prelim}, we give the necessary preliminaries on graphical models and Christoffel polynomials. 
In \Cref{SEC:modifiedCD}, we introduce our coordinate-wise modification $\CDco{n}{\mu}$ of the Christoffel polynomial, and derive some of its properties. 
In \Cref{SEC:ratfunction}, we define rational Christoffel functions (\Cref{DEF:rational-cd}). We show that they satisfy support dichotomy in \Cref{SEC:ratfunction:supdich}. 
Finally, in~\Cref{SEC:ratfunction:regular}, we show that they may be used to estimate the density of the underlying measure after a suitable regularization.

\section{Preliminaries} \label{SEC:prelim}

\subsection{Notations} \label{SEC:notations}
We denote $e_1, e_2, \ldots, e_d$ for the standard basis of $\R^d$. We write $\R[\x] = \R[x_1, x_2, \ldots, x_d]$ for the space of $d$-variate, real polynomials, with the understanding that $x_k = e_k^\top \x$ for $\x \in \R^d$. We write $\N = \{0,1,2,\ldots \}$ for the set of nonnegative integers.
For $\alpha \in \N^d$, we write ${\x^\alpha = x_1^{\alpha_1} \cdot x_2^{\alpha_2} \cdots x_d^{\alpha_d}}$ for the monomials, which form a basis of $\R[\x]$. 
For $n \in \N$, we then write ${\R[\x]_n \subseteq \R[\x]}$ for the subspace spanned by monomials of \emph{total} degree at most $n$, i.e., with $\|\alpha\|_1 := \sum_{k=1}^d \alpha_k \leq n$. 
We write $\Rco \subseteq \R[\x]$ for the subspace spanned by monomials of \emph{coordinate-wise} degree at most $n$, i.e., those with $\|\alpha\|_\infty := {\max_k \alpha_k \leq n}$. 
We have $\mathrm{dim}~\R[\x]_n = s(n, d) := {n + d \choose d}$ and $\mathrm{dim}~\Rco = s(n, d)_\infty := (n+1)^d$. We write $[\x]_n = (\x^\alpha)_{\|\alpha\|_1 \leq n}$ and $[\x]_{\tens{n}} = (\x^{\alpha})_{\|\alpha\|_\infty \leq n}$ for the vectors of monomials of total and coordinate-wise degree at most $n$, respectively.

We write $[d] := \{1,2, \ldots, d\}$. For $A \subseteq [d]$, we write $\R^d_A \subseteq \R^d$ for the subspace spanned by $\{ e_k : k \in A\}$. For $\x \in \R^d$, we write $\x_A = (x_k)_{k \in A} \in \R^d_A$ for the projection of~$\x$ onto~$\R^d_A$. For a subset $\mainset \subseteq \R^d$, we write
$\mainset_A = \{ \x_A : \x \in \mainset \} \subseteq \R^d_A.$ Furthermore, we denote the closure and interior of $\mainset$ by $\clo{\mainset}$, $\inter{\mainset}$, respectively.

\subsection{Conditional independence and graphical models}
\label{SEC:graphmodels}
Let $\mu$ be a probability measure on $\R^d$, assumed to have a density $f$ w.r.t. the Lebesgue measure. For any $A \subseteq [d]$, we denote $\mu_A$ for the marginal of $\mu$ on $\R_A^d$, which has a density $f_A$ w.r.t. the Lebesgue measure on~$\R_A^d$ given by
\begin{equation} \label{EQ:marginal}
f_A(\x_A) := \int f(\x_A, \x_{[d] \setminus A}) d\x_{[d] \setminus A}.
\end{equation}

\begin{definition}
For $A, B, C \subseteq [d]$, we say $\x_A$ is \emph{conditionally independent} of $\x_B$ given $\x_C$ if
\[
    f_C(\x_C) \cdot f(\x) = f_{A \cup C} (\x_{A \cup C}) \cdot f_{B \cup C} (\x_{B \cup C}).
\]
We then write $A \perp_\mu B \,|\, C$. 
\end{definition}

\begin{example}
\label{example:1}
For $d=3$, let $d\mu(\x)=f(\x)d\x$, where $f(\x) = g(x_1,x_2)h(x_2,x_3)$ for some nonnegative measurable functions $g,h$. Then with $A=\{1\}$, $B=\{3\}$ and $C=\{2\}$,
\begin{eqnarray*}
    f_C(x_2)&=&\int g(x_1, x_2)\,h(x_2,x_3)\,dx_1\,dx_3\\
    f_{A \cup C}(\x_{12})&=&\int g(x_1, x_2)\,h(x_2,x_3)\,dx_3\,=\,g(x_1, x_2)\,\int h(x_2,x_3)\,dx_3\\
    f_{B\cup C}(\x_{23})&=&\int g(x_1, x_2)\,h(x_2,x_3)\,dx_1\,=\,h(x_2,x_3)\,\int g(x_1, x_2)\,dx_1\,.
\end{eqnarray*}
\begin{eqnarray*}
f_C(x_2)\,f(\x)&=&f(\x)\,\int g(x_1, x_2)\,h(x_2,x_3)\,dx_1\,dx_3\\
&=&g(x_1, x_2)\,h(x_2,x_3)\,\int g(x_1, x_2)\,h(x_2,x_3)\,dx_1\,dx_3\\
&=&g(x_1, x_2)\,h(x_2,x_3)\,\int g(x_1, x_2)\,dx_1\,\int h(x_2,x_3)\,dx_3\\
&=&f_{A\cup C}(\x_{12})\,f_{B\cup C}(\x_{23})\,,
\end{eqnarray*}
and therefore $A \perp_\mu B \,|\, C$.
\end{example}

Let $G = ([d],~E(G))$ be a simple graph (i.e, without loops or multiple edges between a pair of vertices). For $A, B, C \subseteq [d]$, we say $A$,$B$ are \emph{separated} by $C$ if any path between a vertex $a \in A$ and a vertex $b \in B$ contains a vertex $c \in C$. We say that $G$ is a \emph{graphical model} for $\mu$ if the following holds for all $A, B, C \subseteq [d]$:
\[
    \text{$A, B$ are separated by $C$ in $G$ $\implies A \perp_\mu B \,|\, C$}.
\]
Note that the complete graph on $[d]$ is always a graphical model.
Given a graphical model for $\mu$, our goal is to find a factorization of its density. 
We outline a standard approach to do so using a \emph{chordal completion} of $G$. See, e.g.,~\cite{ProbabilisticGraphicalModels:book, graphicalmodels:book} for a general reference on the use of graphical models in computational statistics, and~\cite{graphmodels} for an overview of relevant graph-theoretical notions.
\begin{definition}
    A graph $G$ is \emph{chordal} if, for any cycle $C$ in~$G$ of length greater than $3$, there is an edge between two vertices in $C$ which is not part of the cycle (called a chord). A \emph{chordal completion} of $G$ is a chordal graph $\overline{G}$ with $V(G) = V(\overline{G})$ and $E(G) \subseteq E(\overline{G})$.
\end{definition} \noindent
Now, choose a chordal completion $\overline{G}$ of the graphical model $G$. 
We denote the maximal cliques of $\overline{G}$ by $\mC = \{C_1, C_2, \ldots, C_\ell\}$. By chordality, it follows that $\ell \leq d$. Furthermore, it follows that these cliques can be arranged in a tree $J$ with vertex set $V(J) = \mathcal{C}$ satisfying the \emph{clique-intersection property}: for any $C \in \mC$ on the (unique) path between $A, B \in \mC$ in $J$, we have ${A \cap B \subseteq C}$. (See \cite[Section 3.1; Theorem 3.2]{graphmodels} for a proof of this fact.) We refer to $J$ as a \emph{junction tree} for~$G$. It should be noted that chordal completions, and junction trees by extension, are not unique in general. In particular, the complete graph is always a chordal completion of $G$ (which leads to a junction tree with one vertex). The upshot is the following.
\begin{proposition} \label{LEM:junctreesep} Let $\mu$ be measure on $\R^d$ with graphical model $G$ and junction tree~$J$. Let $e = \{ A, B\}$ be an edge of $J$. Write $J_A, J_B$ for the two connected components that result from removing $e$ from $J$, 
 and $V_A, V_B \subseteq [d]$ for the variables that occur in $J_A, J_B$, respectively. Then, $\x_{V_A}$ and~$\x_{V_B}$ are independent conditional on $\x_{A \cap B}$.
\end{proposition}
\begin{proof}
The set $A \cap B$ separates $V_A$ from $V_B$ in $G$ (see~\cite[Lemma 4.2]{graphmodels}). Since $G$ is a graphical model for $\mu$, this implies $\x_{V_A}$ and~$\x_{V_B}$ are independent conditional on $\x_{A \cap B}$.
\end{proof}
\medskip\noindent
For instance in Example \ref{example:1} above, $G = (\{1,2,3\}, \{\{1,2\}, \{2,3\}\})$ is a graphical model for $\mu$. It is chordal, with maximal cliques $C_1 = \{1,2\}$ and $C_2 = \{2, 3\}$. Then, $J= (\{C_1, C_2\}, \{ \{C_1, C_2\}\})$ is a junction tree for $G$. Letting $e = \{C_1, C_2\}$ as in the above, we have 
$J_A=\{1,2\}$, $J_B=\{2,3\}$, $V_A=\{x_1, x_2\}$, $V_B=\{x_2,x_3\}$. Hence
$A\cap B=\{2\}$ separates $\{1,2\}$ from $\{2,3\}$ in $G$ and so $\x_{12}$
and $\x_{23}$ are independent conditional on $x_2$.

\medskip\noindent
Recursive application of \Cref{LEM:junctreesep} leads to the following factorization of the density function:
\begin{equation} \label{EQ:juncdecomp0}
\prod_{\{A, B\} \in E(J)} f_{A \cap B}(\x_{A \cap B}) \cdot f(\x) = \prod_{C \in V(J)} f_C(\x_C).
\end{equation}
Assuming positivity of the marginals, this may be rewritten as
\begin{equation} \label{EQ:juncdecomp}
f(\x) = \frac{\prod_{C \in V(J)} f_C(\x_C)}{\prod_{\{A, B\} \in E(J)} f_{A \cap B}(\x_{A \cap B})}.
\end{equation}
Each factor on the right hand side of this equation depends on at most 
\[
\clique(J) := \max_{C \in \mC} |C|
\]
variables. That is, on the size $\clique(\overline{G})$ of the largest clique in $\overline{G}$. 
This number relates to the \emph{treewidth} of the original graphical model, which is defined as $\tau(G) := \min \clique(\overline{G}) - 1$, where the minimum is taken over \emph{all possible} chordal completions $\overline{G}$ of $G$. Treewidth is a well-studied parameter, with deep connections to the computational complexity of various graph problems. 
While chordal completions can be constructed efficiently by adding edges greedily, some cliques of the resulting chordal completion 
may be large. In fact, it is NP-hard to compute the treewidth of a graph, and thus to find a chordal completion that minimizes $\clique(\overline{G})$.

We end this section with an observation that will be useful in later proofs.
One can order the cliques $C_1, C_2, \ldots, C_\ell \in V(J)$ so that they satisfy the \emph{running intersection property (RIP)}\footnote{The running intersection property and clique-intersection property mentioned earlier are sometimes used interchangeably. See~\cite[Theorem 3.5]{graphmodels} for a justification of this ambiguity.}, i.e, so that for any $2 \leq k \leq \ell$, there exists a $j < k$ with
\begin{equation} \label{EQ:RIP}
    \{C_k, C_j\} \in E(J) \text{ and } C_k \cap \big( C_1 \cup C_2 \cup \ldots \cup C_{k-1} \big) \subseteq C_j.
\end{equation}
(In fact, any \emph{reverse topological sorting} of $J$ gives an order with the RIP~\cite[Section 3.3]{graphmodels}.)
We fix one such $j$ for each $k$, and write $S_k := C_k \cap C_j$. This gives the following partition of the variables (where $\dot\cup$ denotes a union of disjoint sets):
\begin{equation} \label{EQ:RIPdisjoint}
    [d] = C_1 \,\dot\cup\, (C_2 \setminus S_2) \,\dot\cup\, \ldots \,\dot\cup\, (C_\ell \setminus S_\ell).
\end{equation}
We can rewrite~\eqref{EQ:juncdecomp0} according to this order as
\begin{equation} \label{EQ:juncdecompRIP0}
    \prod_{k=2}^\ell {f_{S_k}(\x_{S_k})} \cdot f(\x) = \prod_{k=1}^\ell f_{C_k}(\x_{C_k}).
\end{equation} When the marginals are positive, this is equivalent to
\begin{equation} \label{EQ:juncdecompRIP}
    f(\x) = f_{C_1}(\x_{C_1}) \cdot \prod_{k=2}^\ell \frac{f_{C_k}(\x_{C_k})}{f_{S_k}(\x_{S_k})}.
\end{equation}

\begin{example}
In the graphical model depicted in \Cref{FIG:exmp} below, for example, $\x_{125}$ and $\x_{467}$ are independent given~$x_3$. The graph is made chordal by adding an edge $\{\x_4, \x_6\}$. Its junction tree is drawn on the right. Its vertices are the maximal cliques $\{1, 2\}, \{ 2, 5 \}$, $\{2,3\}$, $\{3, 4, 6\}$, and $\{4, 6, 7\}$ of the chordal completion. They satisfy the RIP in that order.
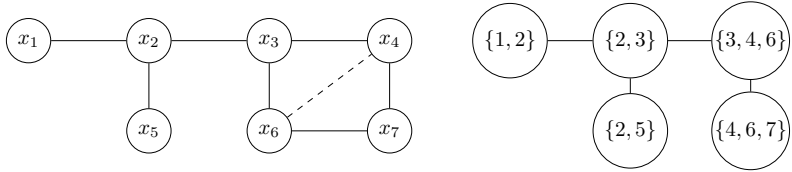
\begin{figure}[H]
\centering
\scalebox{0.8}{
\begin{tikzpicture}
  \node[circle, draw] (A) at (0, 0) {$x_1$};
  \node[circle, draw] (B) at (2, 0) {$x_2$};
  \node[circle, draw] (B2) at (2, -1.5) {$x_5$};
  \node[circle, draw] (C) at (4, 0) {$x_3$};
  \node[circle, draw] (D) at (6, 0) {$x_4$};
  \node[circle, draw] (E) at (4, -1.5) {$x_6$};
  \node[circle, draw] (F) at (6, -1.5) {$x_7$};
  \draw (A) -- (B) -- (C) -- (D) -- (F) -- (E) -- (C);
  \draw (B) -- (B2);
  \draw[dashed] (E) -- (D);
  \begin{scope}[xshift=8cm]
      \node[circle, draw] (A) at (0, 0) {$\{1, 2\}$};
      \node[circle, draw] (B) at (2, 0) {$\{2, 3\}$};
      \node[circle, draw] (B2) at (2, -1.5) {$\{2, 5\}$};
      \node[circle, draw, inner sep=0.2mm] (C) at (4, 0) {$\{3,4,6\}$};
      \node[circle, draw, inner sep=0.2mm] (D) at (4, -1.5) {$\{4,6,7\}$};
      \draw (A) -- (B) -- (C);
      \draw (B) -- (B2);
      \draw (C) -- (D);
  \end{scope}
\end{tikzpicture}
}
\caption{A graphical model with its chordal completion (left) and a corresponding junction tree (right).}
\label{FIG:exmp}
\end{figure}
\end{example}

\subsection{Christoffel polynomials} \label{SEC:CDprelim}
Let $\mainset \subseteq \R^d$ be a compact set with non-empty interior, and let $\mu$ be a probability measure supported on $\mainset$ assumed to have a density $f$ w.r.t. the Lebesgue measure. We have an inner product~$\langle \cdot, \cdot \rangle_\mu$ on the space of polynomials $\R[\x]$ defined by
\begin{equation} \label{EQ:muinner}
    \langle p, q \rangle_\mu := \int p(\x) q(\x) d\mu(\x).
\end{equation}
Any finite-dimensional subspace $V \subseteq \R[\x]$ of dimension $v$ has an orthonormal basis $\{P_1, P_2, \ldots, P_{v}\}$ w.r.t. this inner product, that is, satisfying
\[
    \langle P_i,P_j \rangle_\mu = \delta_{ij} \text{ for all $i,j$}.
\]
This basis gives rise to the so-called \emph{reproducing kernel} of $(V, \langle \cdot, \cdot \rangle_\mu)$, defined as
\[
K^\mu_V(\x, \y) := \sum_{k=1}^{\dim V} P_k(\x) P_k(\y).
\]
This kernel has the property that, for any $p \in V$ and $\x \in \R^n$,
\[
    p(\x) = \langle p, K_V^\mu (\x, \cdot) \rangle_\mu = \int p(\y) K_V^\mu(\x, \y) d \mu(\y).
\]
For the special choice $V = \R[\x]_n$, the kernel $K_{\R[\x]_n}^\mu$ is called the \emph{Christoffel-Darboux kernel} (of degree $n$). The \emph{Christoffel polynomial} is the diagonal of this kernel:
\[
    \CDtot{n}{\mu}(\x) := K_{\R[\x]_n}^\mu(\x, \x) = \sum_{k=1}^{s(n, d)} P_k(\x)^2\,
\]
while the Christoffel \emph{function} is just its reciprocal.
Thus, $\CDtot{n}{\mu}$ is a sum of squares of polynomials of total degree at most $n$. There are several other, equivalent definitions for the Christoffel polynomial. We recall two that will be of use to us later.
\begin{enumerate}
\item \textbf{Via the moment matrix.}
For $\alpha \in \N^d$, recall that we write $\mu_\alpha = \int \x^\alpha d\mu$ for the moments of $\mu$. The (truncated) moment matrix $M_n^\mu \in \mathrm{Mat}(s(n,d))$ of total degree~$n$ for $\mu$ is defined via
\[
    (M_n^\mu \big)_{\alpha, \beta} = \mu_{\alpha + \beta} = \int \x^{\alpha + \beta} d\mu(\x) \quad (\|\alpha\|_1, \,\|\beta\|_1 \leq n).
\]
The Christoffel polynomial $\CDtot{n}{\mu}$ is then given by
\[
    \CDtot{n}{\mu}(\x) = \big([\x]_n\big)^\top (M_n^\mu \big)^{-1} [\x]_n\,,
\]
which is called the ABC Theorem \cite{ABC}.

\item \textbf{Via a variational problem.}
For any $\z \in \R^d$, the Christoffel polynomial can be evaluated at $\z$ by solving the following variational problem:
\begin{equation} \label{EQ:vardef}
    \CDtot{n}{\mu}(\z)^{-1} = \min_{p \in \R[\x]_n} \bigg \{ \int p^2 d\mu : p(\z) = 1 \bigg\}.
\end{equation}
For any  point $\z\in\R^d$, fixed, arbitrary, the optimization problem in the right-hand-side of \eqref{EQ:vardef} 
is a convex quadratic optimization problem that can be solved efficiently even for large dimension $d$. However it only provides the numerical value 
of $\CDtot{n}{\mu}$ at $\z\in\R^d$.
\end{enumerate}
For more details and discussions on the CD kernel and the Christoffel function, the interested reader is referred, e.g., to \cite{dunkl,kroo,Nevai,ABC,Totik}.

\subsection{Asymptotic properties of Christoffel polynomials} \label{SEC:CDsupdich}
Under relatively minor assumptions on $\mainset$ and $\mu$, the Christoffel polynomials exhibit a support dichotomy as $n \to \infty$, captured by the following two theorems from \cite{LPP:2022}. Here, and throughout, we do not consider the behaviour on the boundary $\partial \mainset$ of $\mainset$, which is rather delicate, and lies beyond the scope of this work. See, e.g.,~\cite{LPP:2022, LasserrePauwels2019}.
\begin{theorem}[{\cite[Lemma 4.3.1]{LPP:2022}}]
\label{THM:supdichexp}
    Let $\mu$ be a positive Borel probability measure supported on a compact set $\mainset \subseteq \R^d$. Let $\x \not \in \mainset$, and assume that $\dist(\x, \mainset) \geq \delta$, where $\dist(\cdot, \mainset)$ denotes the Euclidian distance to $\mainset$. Then, for any $n \in \N$, we have
    \[
        \CDtot{n}{\mu}(\x) \geq \frac{1}{8}s(n,d) \cdot 2^{\frac{\delta n}{\delta + \mathrm{diam}(\mainset)}} \cdot n^{-d} \cdot \bigg(\frac{d}{e}\bigg)^d \exp\left(-\frac{d^2}{n}\right).
    \]
\end{theorem}
\begin{theorem}[{\cite[Lemma 4.3.2] {LPP:2022}}, see also~{\cite[Lemma 6.2]{LasserrePauwels2019}}] \label{THM:supdichpoly}
    Let $\mu$ be a probability measure supported on a compact set $\mainset$ with $\mainset = \clo{\inter{\mainset}}$. Assume $\mu$ has a density~$f$ w.r.t. the Lebesgue measure on~$\mainset$, with $f(\x) > c$ for all $\x \in \mathrm{int}(\mainset)$. 
    Let $\x \in \mainset$ and assume that ${\dist(\x, \partial \mainset) \geq \delta}$ for some $\delta > 0$. Then, for any $n \in \N$, we have
    \[
        \CDtot{n}{\mu}(\x) \leq s(n,d) \cdot 2\frac{\lambda(\mainset)}{c \delta^d \omega_d}(1 + d)^3,
    \]
    where $\omega_d$ is the $\mu$-surface area of the $d$-dimensional unit sphere in $\R^{d+1}$.
\end{theorem}
Under some regularity assumptions of $\mainset$ and $\mu$,
more can be said for the asymptotics of $\CDtot{n}{\mu}$ on $\mathrm{int}(\mainset)$ as $n\to\infty$. For instance from \cite[Section 4.4.1] {LPP:2022} one obtains:
\begin{proposition} \label{PROP:denslimit}
Let $\mu_0$ be supported on $\mainset$ and such that, uniformly on compact subsets of $\mathrm{int}(\mainset)$, $\lim_{n\to\infty} \CDtot{n}{\mu_0}(\x)/s(n,d)=W_0(\x)$, where $W_0$ is continuous and positive on $\mathrm{int}(\mainset)$.
If $\mu$ has continuous and positive density 
$g$ w.r.t. $\mu_0$ on $\mathrm{int}(\mainset)$,  then uniformly on compact subsets of $\mathrm{int}(\mainset)$, $\lim_{n\to\infty} \CDtot{n}{\mu}/s(n,d)=W_0(\x)/g$.
In particular, when $\mu_0$ is the equilibrium measure of $\mainset$ 
(with density $\omega_E$ w.r.t. the Lebesgue measure on~$\mainset$) and
$\mu$ has density $f$ w.r.t. the Lebesgue measure on~$\mainset$, then 
\begin{equation}
    \label{EQ:asymptotics-tot}
\lim_{n\to\infty} \frac{1}{s(n,d)}\CDtot{n}{\mu}(\x)\,=\,\frac{\omega_E(\x)}{f(\x)},
\end{equation}
uniformly on compact subsets of $\mathrm{int}(\mainset)$.
\end{proposition}
The equilibrium measure of a compact subset of $\mathbb{C}^d$ originates from potential theory ($d=1$) and pluri-potential theory ($d>1$);
for more details the interested reader is referred, e.g., to  Baran \cite{baran} and Beford and Taylor \cite{bedford}.

\subsection{The regularized Christoffel polynomial} \label{SEC:cdregular}
Importantly, the limit in \eqref{EQ:asymptotics-tot} involves the density $\omega_E$ of the equilibrium measure of $\mainset$ which usually is unknown, except for sets $\mainset$ with special geometry like e.g. the unit box, the unit ball, the simplex, the unit sphere. For instance,
for the interval $[-1,1]$, $\omega_E(x)=1/(\pi\sqrt{1-x^2})$.
As a remedy,
Lasserre~\cite{Lasserre:regularizedCD} introduces a \emph{regularization} of the Christoffel polynomial. For $\varepsilon > 0$, it can be defined via the variational problem as:
\begin{equation}
    \label{EQ:regularized}
    \CDtotreg{n}{\mu}(\z; \epsilon)^{-1} := \min_{p \in \R[\x]_n} \bigg \{ \int p^2 d\mu : \int_{B_\infty(\z, \epsilon)} p(\y)  \frac{d \y}{\epsilon^d} = 1 \bigg\}.
\end{equation}

Here, $B_\infty(\z, \epsilon) = \{ \x \in \R^d : \|\x-\z\|_\infty \leq \frac{1}{2} \epsilon\}$ is the closed $L_\infty$-ball around~$\z$ of diameter $\epsilon$ (whose Lebesgue volume is $\varepsilon^{-d}$).
It turns out that $\CDtotreg{n}{\mu}(\z;\varepsilon)$ can be constructed efficiently again via the moment matrix $M_n^\mu$ by:
\[\CDtotreg{n}{\mu}(\z;\varepsilon)\,=\,([ \x;\varepsilon]_n)^\top(M_n^\mu)^{-1}\,[\x;\varepsilon]_n,\]
where 
\[[\x;\varepsilon]_n\,:=\,\int_{B_\infty(\x,\varepsilon)} [\y]_n\frac{d\y}{\varepsilon^d}\quad\in\,\mathbb{R}[\x,\varepsilon]_n^{s(n, d)},\]
is an explicit vector of polynomials in the variables $(\x,\varepsilon)$, of total degree at most $n$.

The primary motivation for this definition is the following result.
\begin{theorem}[{\cite[Theorem 2.5]{Lasserre:regularizedCD}}] \label{PROP:regdenslimit}
Let $\z \in \R^d$ with $B_\infty(\z, \epsilon) \subseteq \mainset$. Assume that the density $f$ of $\mu$ w.r.t. the Lebesgue measure is continuous and positive on $B_\infty(\z, \epsilon)$. Then there is a $\zeta_\epsilon \in B_\infty(\z, \epsilon)$ with
\[
    \lim_{n \to \infty} \epsilon^{d} \cdot \CDtotreg{n}{\mu}(\z; \epsilon) = \frac{1} {f(\zeta_\epsilon)}.
\]
\end{theorem}
In Theorem \ref{PROP:regdenslimit}, $\varepsilon>0$ is a fixed parameter. However and interestingly, $\CDtotreg{n}{\mu}(\z;\varepsilon)$ can be also viewed as a sum of squares of polynomials in the variables $(\z,\varepsilon)$. 
It tuns out that the asymptotics of $\varepsilon^d \CDtotreg{n}{\mu}(\z;\varepsilon)$ when \emph{both} $n$ and
$1/\varepsilon$ increase, is a delicate issue; for more details the interested reader is referred to  \cite{Lasserre:regularizedCD}.

\section{A coordinate-wise Christoffel polynomial} \label{SEC:modifiedCD}

In this section, we introduce a modification of the Christoffel polynomial of \Cref{SEC:CDprelim} which is adapted to measures with product structure. Our modification is derived from the reproducing kernel on a subspace $V \subseteq \R[\x]$ analogously to the standard definition. However, we now set $V = \Rco$, i.e., $V$ is the space of polynomials of \emph{coordinate-wise} degree at most $n$ (rather than total degree).

\begin{definition} \label{DEF:modifiedCD}
Let $\mu$ be measure on $\R^d$ defining an inner product $\langle \cdot, \cdot \rangle_\mu$ on $\R[\x]$ as in~\eqref{EQ:muinner}. For $n \in \N$, we define the \emph{coordinate-wise Christoffel polynomial} $\CDco{n}{\mu}$ in the following equivalent ways:
\begin{enumerate}
    \item \textbf{Via a reproducing kernel.} Let $K^\mu_{\Rco}$ be the reproducing kernel for $\Rco \subseteq \R[\x]$. Then, \[\CDco{n}{\mu}(\x) := K^\mu_{\Rco}(\x, \x)\,\quad\forall \x\in\R^d\,;\]
    \item \textbf{Via an orthonormal basis.} Let $\{ P_1, \ldots, P_{s(n,d)_\infty}\}$ be an orthonormal basis for $\Rco$ w.r.t. the inner product $\langle \cdot, \cdot \rangle_\mu$.
Then, 
\[
    \CDco{n}{\mu}(\x) := \sum_{k=1}^{s(n,d)_\infty} P_k(\x)^2\,,\quad\forall \x\in\R^d.
\]
\item \textbf{Via a variational problem.} Let $\z \in \R^d$, then $\CDco{n}{\mu}(\z)$ can be evaluated at $\z$ by solving the following variational problem:
\begin{equation} \label{EQ:covardef}
    \CDco{n}{\mu}(\z)^{-1} = \min_{p \in \Rco} \bigg \{ \int p^2 d\mu : p(\z) = 1 \bigg\}\,,\quad\forall\z\in\R^d.
\end{equation}
\item \textbf{Via a moment matrix.} For $\alpha \in \N^d$, recall that we write $\mu_\alpha = \int \x^\alpha d\mu$ for the moments of $\mu$. The (truncated) moment matrix $M_{\tens{n}}^\mu \in \mathrm{Mat}(s(n,d)_\infty)$ of coordinate-wise degree~$n$ for $\mu$ is defined via
\[
    (M_{\tens{n}}^\mu \big)_{\alpha, \beta} = \mu_{\alpha + \beta} = \int \x^{\alpha + \beta} d\mu(\x) \quad (\|\alpha\|_\infty, \,\|\beta\|_\infty \leq n).
    \]
The coordinate-wise Christoffel polynomial $\CDco{n}{\mu}$ is then given by
\begin{equation}
    \label{EQ:ABC}
    \CDco{n}{\mu}(\z) = \big([\z]_{\tens{n}}\big)^\top (M_{\tens{n}}^\mu \big)^{-1} [\z]_{\tens{n}}\,,\quad\forall\z\in\R^d.
\end{equation}
Indeed \eqref{EQ:covardef} reads
\[
\CDco{n}{\mu}(\z)^{-1} = \min_{\mathbf{p} \in \R^{s(n,d)_\infty}} \bigg \{ \mathbf{p}^\top M_{\tens{n}}^\mu\mathbf{p} : \mathbf{p}^\top [\z]_{\tens{n}} = 1 \bigg\},
\]
which is a convex quadratic optimization problem. Applying the standard Karush-Kuhn-Tucker optimality conditions, simple linear algebra yields~\eqref{EQ:ABC}.
\end{enumerate}

\end{definition}

The polynomial $\CDco{n}{\mu}$ is thus a sum of squares of polynomials of coordinate-wise degree at most $n$.
Its coefficients can be computed by inverting a moment matrix, now of size $\mathrm{dim}~\Rco = s(n,d)_\infty = (n+1)^d$. For given $n, d$, it is thus more difficult to compute than the standard Christoffel polynomial. 
Furthermore, while the standard Christoffel polynomial is equivariant under distance-preserving transformations of $\R^d$, our coordinate-wise definition is not. That is, for any isometry $T : \R^n \to \R^n$, we have $\CDtot{n}{\mu \circ T}(\x) = \CDtot{n}{\mu}(T\x)$, but the same is not true for the coordinate-wise Christoffel polynomial. In light of the variational definitions~\eqref{EQ:vardef},~\eqref{EQ:covardef}, the point is that, for any polynomial $p \in \R[\x]$, the composition $p \circ T$ is a polynomial of the same \emph{total} degree, but not necessarily of the same \emph{coordinate-wise} degree. Another way to look at this is that the definition of coordinate-wise degree depends on some explicit choice of generators $x_1, x_2, \ldots, x_n$ of the polynomial ring.

However, if one is interested in exploiting product structure on the measure $\mu$, the coordinate-wise definition appears to be the right one. In particular, it satisfies the decomposition~\eqref{EQ:productdecomp}, whereas the standard Christoffel polynomial does not. 
\begin{proposition} \label{LEM:CDprod}
    Let $\mu$ be a product measure, i.e., with $d\mu(\x) = \prod_{i=1}^d d\mu_i(x_i)$. Then we have
    \[
        \CDco{n}{\mu}(\x) = \prod_{i=1}^d \CDtot{n}{\mu_i}(x_i). 
    \]
\end{proposition}
\begin{proof}
For each $i \in [d]$, let $\{ P_{i, 0}, P_{i, 1}, \ldots, P_{i, n} \}$ be an orthonormal basis for $\R[x_i]_n$ w.r.t. $\langle \cdot, \cdot \rangle_{\mu_i}$, so that $\CDtot{n}{\mu_i}(x_i) = \sum_{j=0}^n P_{i, j}(x_1)^2$. For $\alpha \in \N^d$ with $\|\alpha\|_\infty \leq n$, let $P_\alpha(\x) := \prod_{i=1}^d P_{i, \alpha_i}(x_i)$, which is a polynomial of coordinate-wise degree at most~$n$. Then, for any $\alpha, \beta \in \N^d$ with $\|\alpha\|_\infty, \|\beta\|_\infty \leq n$, we find that
\[
    \int P_\alpha(\x) P_\beta(\x) d\mu(\x) = \prod_{i=1}^d \int P_{i, \alpha_i}(x_i)P_{i, \beta_i}(x_i) d\mu_i(x_i) = \delta_{\alpha\beta}.
\]
That is, the $P_\alpha$ form an orthonormal basis for $\Rco$ w.r.t. $\langle \cdot, \cdot \rangle_\mu$, and so
\begin{align*}
    \CDco{n}{\mu}(\x) &= \sum_{\|\alpha\|_\infty \leq n} P_\alpha(\x)^2 = \sum_{\|\alpha\|_\infty \leq n} \prod_{i=1}^d P_{i, \alpha_i}(x_i)^2 \\
    &= \prod_{i=1}^d \sum_{j = 0}^n P_{i, j}(x_i)^2 = \prod_{i=1}^d \CDtot{n}{\mu_i}(x_i).
\end{align*}
\end{proof}

On the other hand, consider the standard Christoffel polynomial for a measure $d\mu(\x) =d\mu_1(x_1) \cdot d\mu_2(x_2)$ on $\R^2$. Let $\{ P_{1, j} : 0 \leq j \leq n\}$ and $\{ P_{2, j} : 0 \leq j \leq n\}$ be orthonormal bases for $\R[x_1]_n$ and $\R[x_2]_n$ w.r.t. $\langle \cdot, \cdot \rangle_{\mu_1}$, $\langle \cdot, \cdot \rangle_{\mu_2}$ respectively, and assume w.l.o.g. that $\mathrm{deg}(P_{i, j})= j$ for each $i, j$.  As we have seen, the polynomials
$P_{1, i}(x_1)P_{2, j}(x_2)$, $0 \leq i,j \leq n$ are orthonormal w.r.t.~$\langle \cdot, \cdot \rangle_\mu$. It follows that:
\begin{align*}
\CDtot{n}{\mu}(x_1,x_2) &= \sum_{i+j\leq n}P_{1, i}(x_1)^2\,P_{2, j}(x_2)^2 = \CDtot{n}{\mu_1}(x_1)\cdot\left[\sum_{j=0}^n\frac{\CDtot{n-j}{\mu_1}(x_1)}{\CDtot{n}{\mu_1}(x_1)}\cdot P_{2, j}(x_2)^2\right] \\ 
&= \CDtot{n}{\mu_1}(x_1)\cdot\left[\sum_{j=0}^n \theta_j(x_1)\,P_{2,j}(x_2)^2\right],
\end{align*}
where $x\mapsto \theta_j(x)\in (0,1)$ for each $j\leq n$, and $\sum_{j=0}^n\theta_j(x)=1$ for all $x \in \R$.
The term $\sum_{j=0}^n\theta_j(x_1)\,P_{2,j}(x_2)^2$ does not equal $\CDtot{n}{\mu_2}(x_2)$ in general. More results on disintegration properties of the standard Christoffel polynomial can be found in~\cite{cras-dis}.

\subsection{Asymptotics} \label{SEC:modifiedCD:asymptotics}
The asymptotic properties of the coordinate-wise and standard Christoffel polynomial are closely related in light of the following observation.
\begin{remark} \label{REM:CDwedge}
For any $n, d \in \N$, we have that
\[
    \R[\x]_n \subseteq \Rco \subseteq \R[\x]_{nd}.
\]
Using the variational definitions~\eqref{EQ:vardef}, \eqref{EQ:covardef}, this implies that, for any $\x \in \R^{d}$,
\[
    \CDtot{n}{\mu}(\x) \leq \CDco{n}{\mu}(\x) \leq \CDtot{nd}{\mu}(\x).
\]   
\end{remark}
Using \Cref{REM:CDwedge}, the support dichotomy property of the standard Christoffel polynomial captured by \Cref{THM:supdichexp} and \Cref{THM:supdichpoly} translate to the coordinate-wise version readily. For ease of reference, we restate the results here. We do not give explicit constants, rather emphasizing the asymptotic behaviour as $n \to \infty$. We also (slightly) relax the assumptions to make the results more suited to our later applications in \Cref{SEC:ratfunction}.
\begin{corollary} \label{COR:maxdegreeSDexp}
Let $\mu$ be a positive Borel probability measure supported on a compact set $\mainset \subseteq \R^d$. Let $\x \not \in \mainset$ be fixed. As $n \to \infty$, there is an $\alpha > 0$ so that
    \[
        \CDco{n}{\mu}(\x) \geq \Omega\big(\exp(\alpha \cdot n)\big).
    \]
\end{corollary}
\begin{corollary} \label{COR:maxdegreeSDpoly}
Let $\mu$ be a probability measure supported on a compact set with $\mainset = \clo{\inter{\mainset}}$.
Assume $\mu$ has a density $f$ w.r.t. the Lebesgue measure on~$\mainset$ which is continuous and positive on $\inter{\mainset}$. Let $\x \in \inter{\mainset}$ fixed. Then, as $n \to \infty$,
\[
        \CDco{n}{\mu}(\x) \leq O(n^d).
\]
\end{corollary}
\begin{proof}
    Note that we cannot invoke \Cref{THM:supdichpoly} directly, as we do not assume a uniform lower bound on the density $f$ of $\mu$. However, as $\x \in \mathrm{int}(\mainset)$, and $f$ is continuous and positive on $\mathrm{int}(\mainset)$, we may choose a constant $c > 0$ and an open ball $\x \in B \subseteq \mainset$ so that $f(\z) > c$ for all $\z \in B$. 
    Let $\mainset' = \clo{B} \subseteq \mainset$, and let $\mu' = \mu_{|\mainset'} / \mu(\mainset')$ be the (normalized) restriction of $\mu$ to~$\mainset'$. Now, $\mainset'$ and $\mu'$ satisfy the assumptions of \Cref{THM:supdichpoly}.
    Furthermore, it follows from the variational definition~\eqref{EQ:covardef} that $\CDco{n}{\mu}(\x) \leq {\CDco{n}{\mu'}(\x)} \cdot \mu(\mainset')^{-1}$.    
    Using \Cref{REM:CDwedge} and \Cref{THM:supdichpoly}, we conclude that
    \[
        \CDco{n}{\mu}(\x) \leq \frac{\CDco{n}{\mu'}(\x)}{\mu(\mainset')} \leq \frac{\CDtot{nd}{\mu'}(\x)}{\mu(\mainset')} \leq O(n^d).
    \]
\end{proof}
Somewhat surprisingly, it is not clear whether \Cref{PROP:denslimit} generalizes as well.\footnote{In a private communication~\cite{Levenberg:personal}, N. Levenberg suggests that in general the limit of $\CDco{n}{\mu}/(n+1)^d$ as $n\to\infty$, if it exists, is not 
identical to that of $\CDtot{n}{\mu}/s(n, d)$. Moreover, tools to analyze the latter are not available for the former. (Again recall that here $\CDtot{n}{\mu}$ and $\CDco{n}{\mu}$ refer to Christoffel polynomials, i.e., reciprocals of Christoffel functions.)} However we make the following conjecture:
\begin{conjecture} \label{PROP:codenslimit}
Assume that $\mu$ has a positive, continuous density $f$ w.r.t. the Lebesgue measure on $\mainset$. Then we have:
\[
    \lim_{n \to \infty} \frac{1}{s(n,d)_\infty} \cdot \CDco{n}{\mu}(\z) = \frac{\omega_{E, \infty}(\z)} {f(\z)},
\]
where $\omega_{E, \infty} : \R^d \to \R_{\geq 0}$ depends on $\mainset$ (but not on $f$).
\end{conjecture}

\subsection{Regularization} 
Next, as in~\Cref{SEC:cdregular}, one can regularize $\CDco{n}{\mu}$ as follows:
\begin{equation}
   \label{EQ:regularized-new} 
    \CDcoreg{n}{\mu}(\z; \epsilon)^{-1} := \min_{p \in \Rco} \bigg \{ \int p^2 d\mu : \int_{B_\infty(\z, \epsilon)} p(\y)  \frac{d \y}{\epsilon^d} = 1 \bigg\}\,,
\end{equation} 
which is the exact analogue for $\CDco{n}{\mu}$ of \eqref{EQ:regularized} for $\CDtot{n}{\mu}$. 
Again, $\CDcoreg{n}{\mu}(\z;\varepsilon)$ can be constructed efficiently via the moment matrix $M_n^\mu$ by:
\[\CDcoreg{n}{\mu}(\z;\varepsilon)\,=\,([ \x;\varepsilon]_{\tens{n}})^\top (M_n^\mu)^{-1}\,[\x;\varepsilon]_{\tens{n}},\]
\[ [\x;\varepsilon]_{\tens{n}}
\,:=\,\int_{B_\infty(\x,\varepsilon)} ([\y]_{\tens{n}})\,\frac{d\y}{\varepsilon^d}\quad\in\,\mathbb{R}[\x,\varepsilon]_{\tens{n}}^{s(n,d)_{\infty}},\]
is an explicit vector of sums of squares of polynomials in the variables $(\x,\varepsilon)$, of coordinate-wise degree at most $n$.
We then have the following analog of~\Cref{PROP:regdenslimit}.
\begin{theorem} \label{PROP:coregdenslimit}
Let $\varepsilon>0$ be fixed and let $\z \in \R^d$ with $B_\infty(\z, \epsilon) \subseteq \mainset$. Assume that the density $f$ of $\mu$ w.r.t. the Lebesgue measure is continuous and positive on $B_\infty(\z, \epsilon)$. Then there is a $\zeta_\epsilon \in B_\infty(\z, \epsilon)$ with
\[
    \lim_{n \to \infty} \epsilon^{d} \cdot \CDcoreg{n}{\mu}(\z; \epsilon) = \frac{1} {f(\zeta_\epsilon)}.
\]
\end{theorem}
\begin{proof}
Similarly to \Cref{REM:CDwedge}, we have for any $\x \in \R^{d}$ and $\epsilon > 0$ that,
\[
    \CDtotreg{n}{\mu}(\x; \epsilon) \leq \CDcoreg{n}{\mu}(\x; \epsilon) \leq \CDtotreg{nd}{\mu}(\x; \epsilon).
\]
We note that $\big(\epsilon^d \CDtotreg{nd}{\mu}(\x; \epsilon)\big)_{n \in \N}$ is a subsequence of $\big(\epsilon^d \CDtotreg{n}{\mu}(\x; \epsilon)\big)_{n \in \N}$, meaning that
\[
    \lim_{n \to \infty } \epsilon^d \CDtotreg{nd}{\mu}(\x; \epsilon) = \lim_{n \to \infty } \epsilon^d \CDtotreg{n}{\mu}(\x; \epsilon),
\]
whence  $\lim_{n \to \infty } \epsilon^d \CDcoreg{n}{\mu}(\x; \epsilon) = \lim_{n \to \infty } \epsilon^d \CDtotreg{n}{\mu}(\x; \epsilon) = f(\zeta_\epsilon)^{-1}$, by \Cref{PROP:regdenslimit}.
\end{proof}

\subsection{Technical properties.}
We end this section with two technical properties of the coordinate-wise Christoffel polynomial that will be useful later. 

\begin{lemma}
\label{LEM:CDgeq1}
If $\mu$ is a probability measure, then $\CDco{n}{\mu}(\z) \geq 1$ for all $\z \in \R^d$.
\end{lemma}
\begin{proof}
    The constant polynomial $p=1$ is a feasible solution to the variational problem defining $\CDco{n}{\mu}$ for any $n \in \N$. As a result, we find for any $\z \in \R^d$, that
    \[
        \frac{1}{\CDco{n}{\mu}(\z)} \leq \int 1^2 d \mu = 1, \quad \text{ meaning } \CDco{n}{\mu}(\z) \geq 1.
    \]
\end{proof}
\begin{lemma}
\label{LEM:CDfracgeq1}
For any $\emptyset \neq D \subseteq [d]$, we have
\[
\frac{\CDco{n}{\mu}(\z)}{\CDco{n}{\mu_D}(\z_D)} \geq 1 \quad \text{ for all $\z \in \R^d$}.
\]
\end{lemma}
\begin{proof}
    Let $p_D \in \R[\x_D]_{\tens{n}}$ be the minimizer of the variational problem~\eqref{EQ:covardef} defining~$\CDco{n}{\mu}(\z_D)^{-1}$. Then the extension $p(\x) = p_D(\x_D)$ of $p_D$ to $\Rco$ is a feasible solution to the variational problem defining $\CDco{n}{\mu}(\z)^{-1}$, with
    \[
        \int p^2 d\mu = \int p_D^2 d\mu_D.
    \]
    This shows that $\CDco{n}{\mu}(\z)^{-1} \leq \CDco{n}{\mu_D}(\z_D)^{-1}$, and so $\CDco{n}{\mu}(\z) \geq \CDco{n}{\mu_D}(\z_D)$.
\end{proof}

\section{Rational Christoffel functions}
\label{SEC:ratfunction}
Let $\mu$ be a probability measure on $\R^d$, and write $f$ for its density w.r.t. the Lebesgue measure. Assume that $G$ is a graphical model for $\mu$ with junction tree $J$. Recall the decomposition~\eqref{EQ:juncdecomp} of the density function:
\[
f(\x) = \frac{\prod_{C \in V(J)} f_C(\x_C)}{\prod_{\{A, B\} \in E(J)} f_{A \cap B}(\x_{A \cap B})}.
\]
Inspired by this decomposition,
we define the following \emph{rational Christoffel function}:
\begin{definition} \label{DEF:rational-cd}
Let $\mu$ be a probability measure on $\R^d$ with graphical model $G$ and junction tree~$J$. For $\x \in \R^d$, $n \in \N$, we define
\begin{equation}
    \label{EQ:def-rational}
    \x\mapsto \sparseCD{n}{\mu}{G,J}(\x) := \frac{\prod_{C \in V(J)} \CDco{n}{\mu_C}(\x_C)}{\prod_{\{A, B\} \in E(J)} \CDco{n}{\mu_{A \cap B}}(\x_{A \cap B})}.
\end{equation}
\end{definition}
\noindent We can make the following observations immediately:
\begin{enumerate}
\item As emphasized in its notation, the rational function in \eqref{EQ:def-rational} depends on the chosen junction tree $J$ of $G$.
    \item The function $\sparseCD{n}{\mu}{G,J}$ is a rational sum of squares. That is, each of the factors occuring in the numerator and denominator of~\eqref{EQ:def-rational} is a sum of squares of polynomials. In fact, they are Christoffel polynomials of coordinate-wise degree $2n$ associated to marginal distributions $\mu_D$ of $\mu$, with $|D| \leq \clique(J)$ ($= \tau(G) + 1$ if $J$ corresponds to an optimal chordal completion of $G$).
    \item The coefficients of each individual factor can therefore be computed by inverting a matrix of size at most $(n+1)^{\clique(J)}$; namely the principal submatrix of $M^\mu_{\tens{n}}$ indexed by the monomials supported on $D$.    
    In total, we need to perform $|V(J)| + |E(J)| \leq 2d$ such inversions. 
    
    \item  The numerator and denominator of~\eqref{EQ:def-rational} are both strictly positive for all $\x \in \R^d$. The levels sets $\{ \x \in \R^d : \sparseCD{n}{\mu}{G,J}(\x) \leq \gamma \}$ of $\sparseCD{n}{\mu}{G,J}$ are thus semialgebraic, as $\sparseCD{n}{\mu}{G,J}(\x) \leq \gamma$ if, and only if,
    \[
         {\prod_{C \in V(J)} \CDco{n}{\mu_C}(\x_C)} \,\leq\, \gamma \cdot {\prod_{\{A, B\} \in E(J)} \CDco{n}{\mu_{A \cap B}}(\x_{A \cap B})}.
    \]
\end{enumerate}

\begin{example} \label{EXM:compgain}
To illustrate the potential gain in computational complexity of our rational functions w.r.t. standard Christoffel polynomials, consider a $d$-dimensional multivariate measure $\mu$ on $\R^d$ where in its associated graphical model, the cliques are just
$C_i:=\{i,i+1\}$, with associated marginal $\mu_{C_i}$, $i=1,\ldots,d-1$.
Then
\[
\sparseCD{n}{\mu}{G,J}(\x)\,=\,\frac{\prod_{i=1}^{d-1}\CDco{n}{\mu_{C_i}}(x_i,x_{i+1})}{\prod_{i=2}^{d-1}\CDtot{n}{\mu_i}(x_i)}\,,\quad\forall\x\in\R^d\,;\:n\in\N.
\]
Then for instance, to detect whether a point $\x\in\mathrm{supp}(\mu)$ via  a test of the form $\sparseCD{n}{\mu}{G,J}(\x)\leq \gamma$, for some scalar $\gamma$, one only need to compute $d-1$ bivariate Christoffel polynomials $\CDco{n}{\mu_{C_i}}(x_i,x_{i+1})$, and $d-2$ univariate Christoffel polynomials $\CDtot{n}{\mu_{i}}(x_i)$. So we have to invert $d-1$ moment 
matrices of size $O(n^2)$
and $d-2$ moment 
matrices of size $O(n)$, a significant computational gain when compared to inverting a moment matrix of size $O(n^d)$!
\end{example}
As alluded to earlier, the following proposition shows that $\sparseCD{n}{\mu}{G, J}$ is, in some sense, a natural generalization of the modified Christoffel polynomial $\CDco{n}{\mu}$.

\begin{proposition} \label{PROP:agreewhenproduct} Suppose that $d\mu(\x) = \prod_{i=1}^d d\mu_i(x_i)$ is a product measure. Let~$G$ be \emph{any} graphical model for $\mu$, and $J$ any junction tree for $G$.  Then we have
\[
    \sparseCD{n}{\mu}{G,J}(\x) = \prod_{i=1}^d \CDtot{n}{\mu_i}(x_i) = \CDco{n}{\mu}(\x).
\]
\end{proposition}
\begin{proof}
We order the factors of the rational Christoffel function according to \eqref{EQ:juncdecompRIP0}, and then decompose each factor using \Cref{LEM:CDprod}, yielding
\begin{align*}
\sparseCD{n}{\mu}{G,J}(\x)
&= \frac{\prod_{k=1}^\ell \CDco{n}{\mu_{C_k}}(\x_{C_k})}{\prod_{k=1}^\ell \CDco{n}{\mu_{S_k}}(\x_{S_{k}})}
= \frac{\prod_{k=1}^\ell \prod_{i \in C_k} \CDtot{n}{\mu_i}(x_i)}{\prod_{k=1}^\ell \prod_{i \in S_k} \CDtot{n}{\mu_{i}}(x_{i})} \\
&= \prod_{k=1}^\ell ~\prod_{i \in C_k \setminus S_k} \CDtot{n}{\mu_{i}}(x_{i})
= \prod_{i = 1}^d \CDtot{n}{\mu_{i}}(\x_{i}) = \CDco{n}{\mu}(\x),
\end{align*}
where we have used~\eqref{EQ:RIPdisjoint} for the second to last equality.
\end{proof}

\subsection{Support dichotomy} \label{SEC:ratfunction:supdich}
Next, we prove that our rational Christoffel functions satisfy a support dichotomy property. For this, we need the following assumptions on the measure $\mu$. They can be thought of as an extension of the assumptions made in \Cref{SEC:modifiedCD:asymptotics} to all marginals of $\mu$ (compare also to~\cite{LPP:2022, LasserrePauwels2019}).

\begin{assumption} \label{ASSU:measure}
    Let $\mu$ be a probability measure with compact support $\mainset \subseteq \R^d$ and density $f$ w.r.t. the Lebesgue measure. For every $\emptyset \neq A \subseteq [d]$, write $f_A$ for the marginal density of $\mu$ on the variables in $A$ as in~\eqref{EQ:marginal}. We assume that
    \begin{itemize} 
        \item $\mainset_A = \clo{\inter{\mainset_A}}$;
        \item $f_A$ is continuous and positive on $\inter{\mainset_A}$;
        \item $f_A = 0$ outside of $\mainset_A$.
    \end{itemize}
\end{assumption}

\begin{remark} \label{REM:f>0}
    For any $\mainset, \mu$ satisfying~\Cref{ASSU:measure}, and $\emptyset \neq A \subseteq [d]$, we have 
    \[
        \mathrm{int} (\mainset_A) = \{ \x \in \R_A^d : f_A(\x) > 0 \} \setminus \partial \mainset_A.
    \]
\end{remark}
Before we move to the proof of support dichotomy, we wish to first remark that~\Cref{ASSU:measure} comes at a loss of generality w.r.t. the assumptions made in~\Cref{SEC:modifiedCD:asymptotics}, as the following example shows.
\begin{example}
    Consider the `hourglass' $\mainset \subseteq \R^2$, given by
    \[
        \mainset := \{ (x_1, x_2) \in \R^2 : |x_1| \leq |x_2| \leq 1\}.
    \]
    Let $\mu$ be the restriction of the Lebesgue measure to $\mainset$, and let $f \propto \mathbf{1}\{\mathrm{int}(\mainset)\}$ its density. Then, the conditions of~\Cref{ASSU:measure} are satisfied for $A = \{1, 2\}$, but \emph{not} for $A = \{2\}$. Indeed, note that $\mainset_2 = [-1, 1]$, so that $0 \in \mathrm{int}(\mainset_2)$, but $f_2(0) = 0$. Furthermore, note that while $f \geq 1$ is bounded away from zero on~$\mathrm{int}(\mainset)$, the marginal density $f_1$ is arbitrarily close to $0$ near $x_1=1 \in \mathrm{\mainset_1} = [-1, 1]$. So, $f_1$ would not satisfy the conditions of~\Cref{THM:supdichpoly}.
\end{example}

We can now state our main result.
\begin{theorem}\label{THM:sparsesupdich}
Let $\mu$ and $\mainset$ be such that \Cref{ASSU:measure} holds. Let $G$ be a graphical model for $\mu$, with junction tree $J$.
Let $\x \in \R^d$ be fixed. Assume that, for any ${\emptyset \neq A \subseteq [d]}$, $\x_A \not \in \partial \mainset_A$. Then, as $n \to \infty$, there is a constant $\alpha > 0$ such that
\[
    \sparseCD{n}{\mu}{G,J}(\x) =
    \begin{cases} 
        O(n^{d^2}) \quad &\text{if $\x \in \mathrm{int}(\mainset)$},
        \\
        \Omega(\exp(\alpha \cdot n)) \quad &\text{if $\x \not\in \mainset$}.
    \end{cases}
\]
\end{theorem}
As we will see below (cf. \Cref{LEM:supportdecomp}), the proof of \Cref{THM:sparsesupdich} relies crucially on the clique-intersection property of the junction tree $J$. The high-level idea is as follows. We order the factors of the rational Christoffel function as in \eqref{EQ:juncdecompRIP}, i.e.,
\[
    \sparseCD{n}{\mu}{G, J}(\x)
    = \CDco{n}{\mu_{C_{1}}}(\x_{C_1}) \cdot \prod_{k=2}^\ell \frac{\CDco{n}{\mu_{C_k}}(\x_{C_k})}{\CDco{n}{\mu_{S_k}}(\x_{S_k})},
\]
and consider the individual factors
\begin{equation}
\label{EQ:REVIS:Fk}
F_1(\x_{C_1}) := \CDco{n}{\mu_{C_{1}}}(\x_{C_1}), \quad 
F_k(\x_{C_k}) := \frac{\CDco{n}{\mu_{C_k}}(\x_{C_k})}{\CDco{n}{\mu_{S_k}}(\x_{S_k})}, \quad 2 \leq k \leq \ell.
\end{equation}
We then proceed as follows:

\smallskip\noindent
\textbf{Case where $\x \in \mathrm{int}(\mainset)$}: we wish to show that $\sparseCD{n}{\mu}{G, J}(\x) \leq O(n^{d^2})$. To do so, we show upper bound in $O(n^d)$ on each individual factor $F_k$, $1 \leq k \leq \ell$, which follows readily from support dichotomy for (modified) Christoffel polynomials (cf. \Cref{COR:maxdegreeSDpoly}).

\smallskip\noindent
\textbf{Case where $\x \not \in \mainset$}: we need to show that $\sparseCD{n}{\mu}{G, J}(\x) \geq \Omega(\exp(\alpha \cdot n))$. By \Cref{LEM:CDfracgeq1}, we have $F_k(\x) \geq 1$ for all $1 \leq k \leq \ell$, and so it suffices to find a single $\kappa$ such that $F_\kappa(\x) \geq \Omega(\exp(\alpha \cdot n))$. Note that it is possible to find a $\kappa$ with $\x_{C_\kappa} \not \in \mainset_{C_\kappa}$. It follows by support dichotomy for Christoffel polynomials that the \emph{numerator} of $F_\kappa$ is exponentially large in $n$. However, for $\kappa \geq 2$, it could be the case that the \emph{denominator} of $F_\kappa$ is exponentially large as well, meaning we would not be able to conclude anything about $F_\kappa$ itself. It turns out that we can in fact find a $\kappa$ for which the numerator is exponentially large, but the denominator is at most $O(n^d)$. The key to finding such $\kappa$ is the following lemma, which exploits the clique-intersection property of the junction tree.

\begin{lemma} \label{LEM:supportdecomp}
    Let $\mainset, \mu$ be as  in \Cref{ASSU:measure}. Assume that $\mu$ has a graphical model $G$, with junction tree $J$. Let $C_k$, $1 \leq k \leq \ell$ and $S_k$, $2 \leq k \leq \ell$, be as in \eqref{EQ:juncdecompRIP0}. 
    Let $\x \in \R^d$ with $\x_A \not\in \partial \mainset_A$ for all $\emptyset \neq A \subseteq [d]$. If $\x \not\in \mainset$, then either $\x_{C_1} \not \in \mainset_{C_1}$, or there exists a~$k \in \N$ with $2 \leq k \leq \ell$ such that $\x_{C_k} \not \in \mainset_{C_k}$ and $\x_{S_k} \in \mainset_{S_k}$.
\end{lemma}
\begin{proof} 
Recall the decomposition~\eqref{EQ:juncdecompRIP0} of the density function according to the clique-intersection property (aka running intersection property):
\begin{equation} 
    \prod_{k=2}^\ell {f_{S_k}(\x_{S_k})} \cdot f(\x) = \prod_{k=1}^\ell f_{C_k}(\x_{C_k}).
\end{equation}
As $\x \not \in \mainset$, we have $f(\x) = 0$ by \Cref{REM:f>0}, which implies that there is a $k \in [\ell]$ with $f_{C_k}(\x_{C_k}) = 0$. 
Let $\kappa \in [\ell]$ be the smallest index for which ${f_{C_\kappa}(\x_{C_\kappa}) = 0}$. By \Cref{REM:f>0}, this means $\x_{C_\kappa} \not \in \mainset_{C_\kappa}$ (as $\x_{C_\kappa} \not \in \partial\mainset_{C_\kappa}$ by assumption).
If~${\kappa = 1}$, we are thus done immediately. 
So suppose that $\kappa \geq 2$.
By definition, we have ${S_{\kappa} = C_\kappa \cap C_j}$ for some $j < \kappa$ (see \Cref{SEC:graphmodels}). 
By minimality of $\kappa$, we have~$f_{C_{j}}(\x_{C_{j}}) > 0$, which by \Cref{REM:f>0} means that~$\z_{C_{j}} \in \mainset_{C_{j}}$. But that implies that~$\z_{S_{\kappa}} \in \mainset_{S_\kappa}$, as~$S_\kappa \subseteq C_j$. \qedhere
\end{proof}
\begin{proof}[Proof of \Cref{THM:sparsesupdich}]
    We order the factors of the rational Christoffel function as in \eqref{EQ:juncdecompRIP}, and consider the individual factors $F_k$, $1 \leq k \leq \ell$ as in~\eqref{EQ:REVIS:Fk}.
    
    \smallskip \noindent \textbf{Case where $\x \in \mathrm{int}(\mainset)$.} As $\mu_A$ is a probability measure for all $\emptyset \neq A \subseteq [d]$, we may use \Cref{LEM:CDgeq1}, to conclude for all $1 \leq k \leq \ell$ that
    \[
    F_k(\x_{C_k}) \leq \CDco{n}{\mu_{C_k}}(\x_{C_k}).
    \]
    Now, since $\x \in \mathrm{int}(\mainset)$, we have $\x_A \in \mathrm{int}(\mainset_A)$ for all $A \subseteq [d]$. Furthermore, by \Cref{ASSU:measure}, the density $f_{C_k}$ of $\mu_{C_k}$ is positive on $\mathrm{int}(\mainset_{C_k})$ for all $1 \leq k \leq \ell$. We may thus apply~\Cref{COR:maxdegreeSDpoly} to conclude that
    \[
    F_k(\x_{C_k}) \leq {\CDco{n}{\mu_{C_k}}(\z_{C_k})} \leq O(n^d),
    \]
    leading to $\sparseCD{n}{\mu}{G, J}(\x) = \prod_{k=1}^\ell F_k(\x_{C_k}) \leq O(n^{d \cdot \ell}) \leq O(n^{d^2})$.
    
    \smallskip \noindent \textbf{Case where $\z \not\in \mainset$.} 
    Note first that, by \Cref{LEM:CDfracgeq1}, we have for all $1 \leq k \leq \ell$,
    \[
    F_k(\x_{C_k}) = \frac{\CDco{n}{\mu_{C_k}}(\z_{C_k})}{\CDco{n}{\mu_{S_k}}(\z_{S_k})} \geq 1.
    \]
    It is thus sufficient to find a single $k$ such that $F_k(\x) \geq \Omega(\exp(\alpha \cdot n))$, with $\alpha > 0$.
    We use the $k$ obtained in \Cref{LEM:supportdecomp}. It has the property that $\x_{C_k} \not \in \mainset_{C_k}$ but ${\x_{S_k} \in \mainset_{S_k}}$ (if ${k \geq 2}$), which implies ${\x_{S_k} \in \mathrm{int}(\mainset_{S_k})}$ by our assumptions on $\x$. Because of \Cref{ASSU:measure}, we may apply \Cref{COR:maxdegreeSDexp} to find $\alpha' > 0$ such that
    \[
        \CDco{n}{\mu_{C_k}}(\x_{C_k}) \geq \Omega(\exp(\alpha' \cdot n)).
    \]
    And so for $k=1$, we are done. 
    On the other hand, if $k\geq 2$, we use \Cref{COR:maxdegreeSDpoly} to find that
    \[
        \CDco{n}{\mu_{S_k}}(\x_{S_k}) \leq O(n^d).
    \]
    In conclusion, for any $0 < \alpha < \alpha'$, we have that
    \[
    F_k(\x_{C_k}) \geq \frac{\Omega(\exp(\alpha' \cdot n))}{O(n^d)} \geq \Omega(\exp(\alpha \cdot n)).
    \] 
\end{proof}

\subsection{Regularization} \label{SEC:ratfunction:regular}
Following \Cref{SEC:cdregular}, we define a regularized rational Christoffel function, depending on an additional parameter $\epsilon > 0$.
\begin{definition}
For $\x \in \R^d$ and $\epsilon > 0$, we define
\begin{equation} \label{EQ:rat-regularized}
    \sparseCDreg{n}{\mu}{G, J}(\x; \epsilon) := \frac{\prod_{C \in V(J)} \CDcoreg{n}{\mu_C}(\x_C; \epsilon)}{\prod_{\{A, B\} \in E(J)} \CDcoreg{n}{\mu_{A \cap B}}(\x_{A \cap B}; \epsilon)}.
\end{equation} 
\end{definition}
Note that, in order to obtain~\eqref{EQ:rat-regularized}, we have simply replaced each factor in the rational Christoffel function~\eqref{EQ:def-rational} by its regularization~\eqref{EQ:regularized-new}. 
We have the following analog of~\Cref{PROP:coregdenslimit}.
\begin{theorem} \label{THM:sparseconvdensity} Let $\mu, \mainset$ be such that \Cref{ASSU:measure} holds. Let $G$ be a graphical model for $\mu$ with junction tree $J$.
For any $\z \in \R^d$ with $B_\infty(\z, \epsilon) \subseteq \mainset$, we have that
    \[
        \lim_{n \to \infty} \epsilon^{-d} \cdot \sparseCDreg{n}{\mu}{G, J}(\z; \epsilon)^{-1} = \frac{\prod_{C \in V(J)} f_C(\zeta^C)}{\prod_{\{A, B\} \in E(J)} f_{A \cap B}(\zeta^{A \cap B})},
    \]
    where, for each $D \subseteq [d]$ that occurs in the above, $\zeta^D \in B_\infty(\z_D, \epsilon) \subseteq \R^d_D$. 
\end{theorem}
\begin{proof}
    Note that by definition of the junction tree $J$,
    \[
    \sum_{C \in V(J)} |C| \quad - \sum_{\{A, B\} \in E(J)} |A \cap B| = d.
    \]
    As a result, for any $\epsilon > 0$, we have 
    \[
        \epsilon^{-d} \cdot \sparseCDreg{n}{\mu}{G, J}(\z; \epsilon)^{-1} = \frac{\prod_{\{A, B\} \in E(J)} \epsilon^{|A \cap B|} \cdot \CDcoreg{n}{\mu_{A \cap B}}(\z_{A \cap B}; \epsilon)}{\prod_{C \in V(J)} \epsilon^{|C|} \cdot \CDcoreg{n}{\mu_C}(\z_C; \epsilon)}.
    \]
    We conclude by applying \Cref{PROP:coregdenslimit} to each factor of this expression.
\end{proof}
\noindent
The previous suggests that, when $\epsilon$ is small,
\[
\lim_{n \to \infty} \epsilon^{-d} \cdot \sparseCDreg{n}{\mu}{G, J}(\z; \epsilon)^{-1}\approx f(\z).
\]
Indeed, this can be made precise as follows.
\begin{corollary} \label{COR:sparsedensityeps} Let $\mu, \mainset$ be such that \Cref{ASSU:measure} holds. Assume in addition that the marginal densities $f_A$ are locally Lipschitz on $\inter{\mainset_A}$ for all $
\emptyset \neq A \subseteq [d]$.  Let $G$ be a graphical model for $\mu$ with junction tree $J$.
For any $\z \in \mathrm{int}(\mainset)$, we then have
\[
\lim_{n \to \infty} \epsilon^{-d} \cdot \sparseCDreg{n}{\mu}{G, J}(\z; \epsilon)^{-1} = f(\z) + O(\epsilon), \quad \text{ as } \epsilon \to 0.
\]
\end{corollary}
\begin{proof}
Let $\eta > 0$ such that $B = \clo{B_\infty(\z, \eta)} \subseteq \mathrm{int}(\mainset)$. For all $D \subseteq [d]$, the marginal density $f_D$ is continuous and positive on $\mathrm{int}(\mainset)_D$, and so there exists an $m > 0$ with $f_D(\x_D) > m$ for all $\x \in B$. By the additional assumption that $f_D$ is locally Lipschitz, there is an $L > 0$ such that $|f_D(\z_D) - f_D(\y_D)| \leq L \cdot \|\z - \y\|_\infty$ for all $\y \in B$. (Here, $m, L$ do not depend on $D$). For $0 < \epsilon < \eta$, and $\zeta^D$ as in \Cref{THM:sparseconvdensity}, we find that
    \[
    |f_D(\zeta_D) - f_D(\z)| \leq L \cdot \|\zeta_D - \z\|_\infty \leq L \cdot \epsilon.
    \]
Now, as $\epsilon \to 0$, we find that
\begin{align*}
\frac{\prod_{C \in V(J)} f_C(\zeta^C)}{\prod_{\{A, B\} \in E(J)} f_{A \cap B}(\zeta^{A \cap B})} &= \frac{\prod_{C \in V(J)} \big(f_C(\z_C) + O(\epsilon)\big)}{\prod_{\{A, B\} \in E(J)} \big(f_{A \cap B}(\z_{A \cap B}) + O(\epsilon)\big)} \\
&= \frac{\prod_{C \in V(J)} f_C(\z_C) + O(\epsilon) }{\prod_{\{A, B\} \in E(J)} f_{A \cap B}(\z_{A \cap B}) + O(\epsilon)} \\
&= \frac{\prod_{C \in V(J)} f_C(\z_C)}{\prod_{\{A, B\} \in E(J)} f_{A \cap B}(\z_{A \cap B})} + O(\epsilon) = f(\z) + O(\epsilon). 
\end{align*}
\end{proof}

\section{Discussion}

\subsection*{The empirical setting}
In applications, one usually does not have access to the exact moments of $\mu$. Rather, one has access to a finite set of samples ${X_1, \ldots, X_N \in \R^d}$ drawn independently according to $\mu$, which allow to compute \emph{empirical moments} 
\[
    \widetilde{\mu}_\alpha = \frac{1}{N}\sum_{k=1}^N \big(X_k\big)^\alpha \approx \int_\mainset \x^\alpha d\mu.
\]
As $N$ grows, one expects that the \emph{empirical} Christoffel polynomial $\CDtot{n}{\widetilde\mu}$ corresponding to these moments approximates the (true) Christoffel polynomial associated to~$\mu$ rather well. In~\cite[Theorem 3.9]{LasserrePauwels2019}, this is made precise for certain settings (see also~\cite[Chapter 6]{LPP:2022}). We note that the rational Christoffel functions defined in this work can be adapted to the empirical setting analogously. It is an interesting question for future research if, and how, the results of~\cite{LasserrePauwels2019} would generalize.

\subsection*{Coordinate-wise vs. standard Christoffel polynomials.}
Instead of the rational function $\sparseCDshort{n}{\mu}{G,J}$ defined in~\eqref{EQ:def-rational}, one could also consider the rational function
\begin{equation}
    \label{EQ:rational-intro-classic}
\x\mapsto \Psi^{\mu}_{n}(\x)\, := \frac{\prod_{C \in V(J)} \CDtot{n}{\mu_C}(\x_C)}{\prod_{\{A, B\} \in E(J)} \CDtot{n}{\mu_{A \cap B}}(\x_{A \cap B})},
\end{equation}
based on products of standard (rather than modified) Christoffel polynomials, and
its regularization~$\widetilde{\Psi}^{\mu}_{n}$ as well.
It turns out that (i)
$\Psi^{\mu}_n$ shares the same support dichotomy property as~$\sparseCDshort{n}{\mu}{G,J}$, and, (ii) via the asymptotics of its regularized version, as $n$ increases one may also use it to approximate the density $f$ of $\mu$.

In theory, $\Psi^{\mu}_n$ is an attractive alternative to $\sparseCDshort{n}{\mu}{G,J}$ as its factors are somewhat easier to compute. Our main motivation to consider $\sparseCDshort{n}{\mu}{G,J}$ in this work is the fact that, as mentioned, whenever $\mu$ is a product measure, we have $\sparseCDshort{n}{\mu}{G,J} = \CDco{n}{\mu}$ regardless of the choice of graphical model for $\mu$, and we have seen that depending on this choice, $\Psi^{\mu}_n$ might not be equal to $\CDtot{n}{\mu}$ even in that very special setting.

On the other hand, the asymptotics of appropriately scaled $\Psi_n$ as $n\to\infty$ may agree 
with the asymptotics of $\Lambda^\mu_n/s(n,d)$ when it exists (under some regularity assumptions). 
For illustration, consider the simple example where $\mu=h(x_1,x_2)g(x_2,x_3)dx_1dx_2dx_3$ on a compact set $\mainset\subset\R^3$, hence with $C_1=\{1,2\}$ and $C_2=\{2,3\}$. Assume that under some regularity assumptions,
\begin{eqnarray*}
\lim_{n\to\infty}\frac{\Lambda^{\mu_{C_1}}_n(x_1,x_2)}{s(n,2)}&=&\omega^{C_1}_E(x_1,x_2)/f_{C_1}(x_1,x_2)\\
\lim_{n\to\infty}\frac{\Lambda^{\mu_{C_2}}_n(x_2,x_3)}{s(n,2)}&=&\omega^{C_2}_E(x_2,x_3)/f_{C_2}(x_2,x_3)\\
\lim_{n\to\infty}\frac{\Lambda^{\mu_{C_1\cap C_2}}_n(x_2)}{s(n,1)}&=&
\omega^{C_1\cap C_2}_E(x_2)/f_{C_1\cap C_2}(x_2)\,.
\end{eqnarray*}
Then
\[
    \lim_{n\to\infty}\frac{s(n,1)}{s(n,2)^2}\,\Psi^\mu_n(x_1,x_2,x_3)\,=\,
\frac{\omega^{C_1}(x_1,x_2)\cdot\omega^{C_2}(x_2,x_3)}{\omega^{C_1\cap C_2}(x_2)}\cdot\frac{f_{C_1\cap C_2}(x_2)}{f_{C_1}(x_1,x_2)\,f_{C_2}(x_2,x_3)}\]
\begin{equation}
\label{asymp-1}
=\,\frac{\omega^{C_1}(x_1,x_2)\cdot\omega^{C_2}(x_2,x_3)}{\omega^{C_1\cap C_2}(x_2)}
\cdot\frac{1}{h(x_1,x_2)\cdot g(x_2,x_3)}\,,
\end{equation}
to compare with
\begin{equation}
    \label{asymp-2}
\lim_{n\to\infty}\frac{\Lambda^\mu_n(x_1,x_2,x_3)}{s(n,3)}
\,=\,\omega_E(x_1,x_2,x_3)\cdot \frac{1}{h(x_1,x_2)\cdot g(x_2,x_3)}\,.\end{equation}
As $\lim_{n\to\infty}s(n,1)/s(n,2)^2=\frac{3}{2}\lim_{n\to\infty}s(n,3)$, then
if 
\[\omega_E(x_1,x_2,x_3)\,=\,\frac{\omega^{C_1}(x_1,x_2)\cdot\omega^{C_2}(x_2,x_3)}{\omega^{C_1\cap C_2}(x_2)}\]
(which is true e.g. when $\mainset=[a,b]\times [c,d]\times [d,e]$)
then both limits in \eqref{asymp-1} and \eqref{asymp-2} are equal up to a constant factor, namely the ratio between the equilibrium measure $\omega_E(x_1,x_2,x_3)$ of $\mainset$ and the density of $\mu$. So both 
$\Gamma^\mu_n$ and $\Psi^\mu_n$ have some pros and cons.\\

We hope that the results presented in this paper will elicit some curiosity and further investigation from the research community in approximation theory on the relation between the Christoffel polynomials $\CDtot{n}{\mu}$ and $\CDco{n}{\mu}$.
In particular, so far it is not clear at all what the asymptotic behavior of $\CDco{n}{\mu}/s(n,d)_\infty$ looks like,
whereas under some regularity properties of $\mainset$, $\mu$, the asymptotics of $\CDtot{n}{\mu}/s(n,d)$ are quite well understood (cf. \Cref{PROP:denslimit}).

\subsection*{Acknowledgments}
We thank the anonymous referees for their valuable comments and suggestions.

\bibliographystyle{plain}
\bibliography{sparseCD}

\end{document}